\documentclass[11pt]{article}
\usepackage[utf8]{inputenc}
\usepackage{graphicx} 
\usepackage{amsmath, amsthm, amssymb} 
\usepackage{hyperref} 
\usepackage{multicol}
\usepackage{multirow}
\usepackage{caption}
\usepackage{subcaption}
\usepackage[labelfont=bf]{caption}
\usepackage{pdflscape}
\usepackage{xcolor}

\usepackage{algorithm}
\usepackage[noend]{algpseudocode}

\newcommand{\myalgsheader}[0]

\algnewcommand{\IIf}[1]{\State\algorithmicif\ #1\ \algorithmicthen}
\algnewcommand{\EndIIf}{\unskip\ \algorithmicend\ \algorithmicif}
\algnewcommand{\IElse}[1]{\State\algorithmicelse\ #1\ }
\algnewcommand{\IfThenElse}[3]{
  \State \algorithmicif\ #1\ \algorithmicthen\ #2\ \algorithmicelse\ #3}

\begin{document}
\title{Balancing Inexactness in Mixed Precision Matrix Computations}

\author{Erin Claire Carson\footnote{Faculty of Mathematics and Physics, Charles University, carson@karlin.mff.cuni.cz. \\The author is supported by the European Union (ERC, inEXASCALE, 101075632). Views and opinions expressed are those of the authors only and do not necessarily reflect those of the European Union or the European Research Council. Neither the European Union nor the granting authority can be held responsible for them. The author additionally acknowledges support from the Charles University Research Centre program No. UNCE/24/SCI/005.}}
\date{}
\maketitle

\paragraph{Abstract.} Support for arithmetic in multiple precisions and number formats is becoming increasingly common in emerging high-performance architectures. From a computational scientist's perspective, our goal is to determine how and where we can safely exploit mixed precision computation in our codes to improve performance.
One case where the use of low precision is natural, common in computational science, is when there are already other significant sources of ``inexactness'' present, e.g., discretization error, measurement error, or algorithmic approximation error. In such instances, analyzing the interaction of these different sources of inexactness can give insight into how the precisions of various computations should be chosen in order to ``balance'' the errors, potentially improving performance without a noticeable decrease in accuracy. We present a few recent examples of this approach which demonstrate the potential for the use of mixed precision in numerical linear algebra and matrix computations.

\section{Introduction}
\label{sec:1}

The scale of available computational power is constantly growing. As of 2022, we have entered the ``exascale era'', meaning that we have computers capable of performing $10^{18}$ (a billion billion) double precision floating point operations per second. This level of computing power offers massive opportunity for scientific simulation and discovery, but we still have a significant challenge ahead, in that we need to actually design algorithms and software for real-world applications that can take advantage of this powerful underlying hardware. 

The focus of this manuscript is on basic numerical linear algebra and matrix computations, which form the core of a number of critical large-scale applications, including computational fluid dynamics, climate and weather modeling, image analysis, AI, etc. A critical point is that in all these applications, the matrix computations that we perform are necessarily \emph{inexact} matrix computations. 

There are many different sources of inexactness that arise in the computational science process. We typically start with a given phenomenon to study and we write down a mathematical model that describes it. We of course can't capture the full reality using a model, so we incur some \emph{modeling error}. In order to actually solve our model equations and perform a simulation, we incur a \emph{discretization error} and/or a \emph{linearization error}. We eventually arrive at a linear algebra problem, and must pick an algorithm to solve it. This involves what we will call \emph{algorithmic approximation errors}. For example, we might approximate our operator or use randomization to make the problem computationally feasible, or we might pick a certain stopping criterion for an iterative method based on the needs of the application. Finally, when we go to run this algorithm on real hardware, we necessarily incur \emph{rounding errors} due to the use of finite precision. 

These errors are intimately connected. For example, in order to determine an appropriate stopping criterion for an iterative method, one must have some understanding of its limits in finite precision, and the particular discretization used will result in different matrix properties, which can change the convergence behavior of an iterative solver. Unfortunately, these different sources of error are almost always studied in isolation. Further, the rounding error is frequently ignored completely. 

Not accounting for rounding error may have made sense for some select applications in the past when all computations were run in a uniform double (64-bit) precision. On today's high-performance compute hardware, however, which features a variety of bit-width formats,  ignoring the rounding error and how it interacts with other sources of inexactness is not only potentially hazardous, but is also a lost opportunity for improving the performance of scientific codes by exploiting the capabilities of today's hardware. 

\section{Floating Point Numbers and Mixed Precision Hardware}
\label{sec:2}

In computer hardware as we know it, we only have a finite number of bits to store and compute with numbers. The industry standard are the IEEE floating point numbers. A base-2 floating point number can be represented as 
\[
(-1)^\text{sign} \times 2^{\text{(exponent-offset)}} \times 1.\text{significand}.
\]
There is one sign bit, some number of exponent bits, and some number of bits for the significand (also called mantissa or fraction), depending on the particular format (which also defines the constant offset). 

The number of bits assigned to the exponent will determine the range of representable numbers, and the number of bits assigned to the significand will affect the unit roundoff $u$, which is the precision with which numbers are stored.
For example, for IEEE 754 double precision (64 bits), we have 11 exponent bits with an offset of 1023, and 52 explicitly stored significand bits (plus the implicit leading one). This means that the largest representable number is $2^{(2^{\#\text{exponent bits}-1)}-1023} = 2^{1024}\approx 2\cdot 10^{308}$,
the smallest (normal) representable number is
$2^{1-1023} = 2^{-1022}\approx 2\cdot 10^{-308}$,
and the unit roundoff is
$2^{-\#\text{significand bits}}=2^{-53}\approx 10^{-16}$.

Whenever we perform a computation with two floating point numbers, the result is rounded to a floating point number, and this incurs a relative error that is bounded by the unit roundoff, i.e., for floating point numbers $x$ and $y$, 
\[
\text{fl}(x \text{ op } y) = (x \text{ op } y)(1+\delta), \quad |\delta|\leq u, \quad \text{op}=+,-,*,/,
\]
where $\text{fl}()$ denotes the result of a computation performed in floating point arithmetic. 

On modern high-performance computer hardware, we have many more precisions to work with than just double. For example, on the NVIDIA H100 GPU (featured in, for example, the JUPITER exascale computer), we have many different formats, not all of which are from the IEEE 754 standard. These formats along with their properties are displayed in Table \ref{tab:formats}. From the table we can see that of course, with fewer bits, we generally have less to work with in both the exponent and significand fields, and thus we generally end up with a smaller range of representable numbers and a greater unit roundoff. From the final column, however, we can see that the potential performance gains from using lower precision are massive. 

The performance advantages of low precision have been highlighted through the relatively recent addition of the HPL-MxP benchmark \cite{HPLMxP}, a complement to the HPL benchmark in the TOP500 ranking \cite{TOP500_homepage}. As of writing, the current fastest computer in the world is El Capitan at Lawrence Livermore National Laboratory, which achieves $1.8$ exaflops/s on the HPL benchmark, which is doing Gaussian elimination with partial pivoting to solve a dense linear system to double precision accuracy. The HPL-MxP benchmark also solves a dense linear system to double precision accuracy, but allows the use of \emph{mixed precision} iterative refinement to do so. On this benchmark, El Capitan reaches $16.7$ effective exaflops/s, demonstrating the power that comes from exploiting the available low precision hardware. 

Indeed, this is a growing trend. The number of systems in the TOP500 ranking that feature mixed precision accelerators is only increasing over time, and further, on the top machines, the vast majority of floating point performance comes from accelerators; see, e.g., \cite{chalmers2023optimizing}. If we want to make effective use of these computers, we must thus redesign our algorithms so that they can make use of the mixed precision hardware available on GPUs and other specialized accelerators.

\begin{table}[]
\centering
\caption{\centering Floating point number formats on the NVIDIA H100 GPU \cite{NVIDIAH100}.}
\label{tab:formats}
\begin{tabular}{c|c|c|c|c|}
\cline{2-5}
                               & size(bits) & range          & $u$               & Tflops performance (NVIDIA H100 TC) \\ \hline
\multicolumn{1}{|c|}{fp64}     & 64         & $10^{\pm 308}$ & $1\cdot 10^{-16}$ & 67                                    \\ \hline
\multicolumn{1}{|c|}{fp32}     & 32         & $10^{\pm 38}$  & $6\cdot 10^{-8}$  & 989                                   \\ \hline
\multicolumn{1}{|c|}{tf32}     & 19         & $10^{\pm 38}$  & $5\cdot 10^{-4}$  & 989                                   \\ \hline
\multicolumn{1}{|c|}{fp16}     & 16         & $10^{\pm 5}$   & $5\cdot 10^{-4}$  & 1979                                  \\ \hline
\multicolumn{1}{|c|}{bf16}     & 16         & $10^{\pm 38}$  & $4\cdot 10^{-3}$  & 1979                                  \\ \hline
\multicolumn{1}{|c|}{fp8-e5m2} & 8          & $10^{\pm 5}$   & $1\cdot 10^{-1}$  & 3958                                  \\ \hline
\multicolumn{1}{|c|}{fp8-e4m3} & 8          & $10^{\pm 2}$   & $6\cdot 10^{-2}$  & 3958                                  \\ \hline
\end{tabular}
\end{table}

\section{The Challenges of Low Precision}
\label{subsec:2}

While it is clear we need to use low precision if we want to get anywhere close to peak performance on today's machines, this is not something we can do blindly.

\paragraph{Reduced Numerical Range.}
The first difficulty is that, as seen in Table \ref{tab:formats}, lower precisions have a smaller range of representable numbers. In floating point error analyses, it is often assumed that no overflow or underflow occurs during a computation. When using low precision, however, this assumption may not be valid. Overflow will usually cause a computation to fail entirely; underflow may be innocuous in some cases, but in others, it can cause us to lose important numerical properties, such as nonsingularity or positive definiteness. 

The reduced numerical range of low precision brings challenges on many levels. From a software engineering perspective, one must carefully implement code to catch and resolve these errors; see, e.g., \cite{scott2025developing}. From an analysis and algorithm development perspective, we would like to identify properties of the problem and precision which are likely to result in underflow or overflow, and develop techniques to mitigate this. One example is the sophisticated scaling and shifting approach developed in \cite{higham2019squeezing} to ``squeeze'' a matrix into lower precision. 

\paragraph{Invalid Error Bounds.}

We also must consider whether existing analysis of stability and accuracy for our algorithms still apply when we use finite precision. In many cases, bounds on the backward or forward errors for a particular algorithm will usually hold only under some assumption like $nu<1$, where $n$ is the problem dimension. Using lower precision means that the unit roundoff $u$ is larger, which means that error bounds only hold for much smaller problems. For example, in half precision where $u\approx 10^{-4}$, we have no guarantee of stability for problems with dimension $10^4$ or larger (a relatively small problem in the world of exascale). 

The question arises as to whether this is an actual limitation of the algorithm or the arithmetic system, or whether such conditions are overly restrictive as a result of the worst-case fashion of standard error analyses. For example, it has been shown that, using probabilistic approaches, one can obtain bounds in which $n$ is replaced by $\sqrt{n}$ (or even by 1 in some special cases) \cite{higham2019new}. Further details and a discussion of other techniques are described in the blog post of Higham \cite{Higham2021ExtremeScaleLowPrecisions}.

\paragraph{Reduced Precision.}

Every time we store a number in lower precision, we lose accuracy, and every time we perform a floating point operation in lower precision, we lose accuracy. As mentioned in the introduction, in many instances in scientific computing, it is assumed that the model error, discretization error, etc., dominates the rounding error, and thus the rounding error is ignored in many analyses. When our rounding errors are on the order of $10^{-4}$ instead of $10^{-16}$, this may no longer be a valid assumption. It is thus now more important than ever to carefully analyze these multiple sources of error \emph{together}, and ultimately, to \emph{balance} the errors. We will give three examples of this in Section \ref{sec:balancing}.

\section{A (Very) Brief History of Mixed Precision Iterative Refinement}

Before jumping into our examples of balancing rounding error with other errors, we will give a (very) brief overview of some work on mixed precision iterative refinement. For a more complete history and further references, we direct the reader to Chapter 3 of the Ph.D. thesis of Vieubl\'{e} \cite{Vieuble2022PhD}. 

Iterative refinement, shown in Algorithm \ref{alg:ir}, is an algorithm for iteratively improving the solution to linear systems $Ax=b$. The initial approximate solution in line \ref{line1} is typically obtained via an LU factorization, with the computed LU factors being reused for solving for the correction $d_i$ in line \ref{line4}.

\begin{algorithm}
\caption{Iterative Refinement}\label{alg:ir}
\begin{algorithmic}[1]
\Require nonsingular $A\in\mathbb{R}^{n\times n}$, $b\in \mathbb{R}^n$, integer maxit
\Ensure Approximate solution $x_{i+1}$ to $Ax=b$ 
\State Solve $Ax_0=b$. \label{line1}
\For{$i=$ 0:maxit}
    \State $r_i = b - Ax_i$ \label{line3}
    \State Solve $Ad_i = r_i$. \label{line4}
    \State $x_{i+1}=x_i+d_i$
\EndFor
\end{algorithmic}
\end{algorithm}

Wilkinson and his colleagues were using this approach as early as 1948 for solving linear systems on the Automatic Computing Engine \cite{Wilkinson1948}, and Wilkinson wrote down a finite precision (fixed point) analysis of this algorithm in his 1963 book \cite{Wilkinson1963Rounding}. The traditional approach used by Wilkinson used two precisions: a working precision $u$ for storage and all computations except
 the residual computation in line \ref{line3}, which was done in precision $u^2$ (double the working precision; here and in the remainder of the paper, we use ``precision $u$'' as shorthand for ``precision with unit roundoff $u$''. This use of mixed precision was quite natural, as the hardware at the time computed an exact scalar product in precision $u^2$ of two numbers in precision $u$ at no extra cost. Wilkinson's  analysis showed that if $3nu\kappa_\infty(A)<1$, where $\kappa_\infty(A)$ is the infinity-norm condition number of $A$, then the limiting relative forward and backward errors, i.e., $\Vert x_i-x\Vert_\infty/\Vert x\Vert_\infty$ and $\Vert b-Ax_i\Vert_\infty/(\Vert A \Vert_\infty \Vert x_i \Vert_\infty \Vert b \Vert_\infty)$, are on the order $O(u)$. 

In \cite{carson2017new}, it was shown that by changing the solver in line \ref{line4}, more ill-conditioned systems can be solved. In particular, if we instead use preconditioned GMRES, where the preconditioner consists of the computed LU factors, the condition for convergence now contains a $\kappa_\infty(U^{-1}L^{-1}A)$ instead of $\kappa_\infty(A)$, which can be substantially smaller, even if the LU factors are computed inexactly. This approach is called GMRES-based iterative refinement or simply GMRES-IR. 

Motivated by emerging hardware that implemented half precision formats, \cite{carson2018accelerating} combined Wilkinson's approach with low-precision factorization variants (see, e.g., \cite{Langou2006Exploiting32bit}) in which the LU factorization (the most expensive part of the algorithm) is computed in half the working precision, and the rest of the computations are carried out in precision $u$. This resulted in bounds for a general \emph{three-precision} variant of iterative refinement, which uses precision $u_f$ for the factorization, $u_r$ for the residual computation, and $u$ for other computations, with $u_r \leq u \leq u_f$. The analysis also contains an ``effective solve precision''  which allows for a general solver in line \ref{line4}.

Iterative refinement in general represents one case where it is natural to use mixed precision, that is, when the algorithm itself contains some type of ``self-correction'' mechanism, or some type of ``inner-outer'' solve scheme. This is a common occurrence in numerical algorithms. 
\section{Balancing Sources of Inexactness}
\label{sec:balancing}

Another case where it is natural to use low precision in parts of a computation is when there exist other significant sources of error in the algorithm, for example, the use of low-rank approximations, or coarse approximations of the domain. Again, while most existing analysis of these types of errors typically ignore finite precision errors, if we want to safely and effectively exploit low precision computations, we must look closely at what happens when we combine these types of analyses to determine how these errors interact. Our goal will ultimately be to balance these errors, and we can accomplish this by producing analyses which tell us how large the finite precision error can be relative to the other approximations we make so that it does not dominate.  

Here we will present three examples of this balancing, dealing with sparsified matrices, randomized algorithms, and hierarchical matrix approximations. All numerical experiments were performed using MATLAB 2025a.

\subsection{Example 1: Mixed Precision Sparse Approximate Inverse Preconditioners} %

Sparse approximate inverses are commonly-used algebraic preconditioners for Krylov subspace methods. For a thorough overview, see \cite{benzi2002preconditioning}.
The goal is to construct a sparse matrix $M$ that is an approximation of $A^{-1}$. In the classical approach, due to Grote and Huckle \cite{grote1997parallel}, shown in Algorithm \ref{alg:spai}, each column of $M$ can be constructed independently. While this is in theory highly parallelizable, construction and memory requirements can still be costly for large scale problems; see, e.g., \cite{he2020efficient}. We thus aim to improve the performance and memory requirements by using low precision to construct this sparse approximate inverse if possible. 

\begin{algorithm}
\caption{Sparse Approximate Inverse Construction}\label{alg:spai}
\begin{algorithmic}[1]
\Require nonsingular $A\in\mathbb{R}^{n\times n}$, initial sparsity structure $\mathcal{J}$, tolerance $\tau$
\Ensure Sparse approximate inverse $M\approx A^{-1}$ 
\For{each column $k$}
    \State Compute QR factorization of submatrix of $A$ defined by $\mathcal{J}$. \label{spai:line2}
    \State Use QR to solve $\min_{m_k}\Vert e_k-Am_k\Vert_2$ where $e_k$ is the $k$th column of the identity.
    \State {\bf{if} $\Vert r_k\Vert_2 = \Vert e_k - Am_k\Vert_2 \leq \tau$ {\textbf{then}} $\text{break;}$}
    \State {\textbf{else} Add select nonzeros to $\mathcal{J}$, repeat from line \ref{spai:line2}.}
\EndFor
\end{algorithmic}
\end{algorithm}

In \cite{carson2023mixed}, it is shown that we can execute the entire computation in some precision $u_s$ and still reach a solution with $\Vert \hat{r}_k \Vert_2 \leq \tau$ (where $\hat{\cdot}$ denotes a computed quantity) if
\[
n^3 u_s \Vert |e_k| + |A| |\hat{m}_k| \Vert_2 \leq \tau.
\]
In other words, the problem must not be so ill-conditioned with respect to $u_s$ that we incur an error greater than $\tau$ just computing the residual if we want the stopping criterion to be satisfiable. 

We can turn this into a looser but perhaps more useful a priori bound, and say that the sparse approximate inverse can be constructed in precision $u_s$ as long as 
\[
u_s \text{cond}_2(A) \lesssim \tau, 
\]
where $\text{cond}_2(A) = \Vert |A^{-1}||A|\Vert_2$.\footnote{Here we use $\lesssim$ to indicate that this is more of a heuristic than a rigorous bound, since we have dropped dimensional constants (which are usually overestimates in practice anyway).}
This tells us how our errors, coming from the stopping criterion in sparse approximate inverse construction, and the rounding error from the algorithm, must be \emph{balanced}. For a given matrix $A$ and a desired value of $\tau$, we have a limit on how large we can make $u_s$. 

The results in \cite{carson2023mixed} also show that this sparse approximate inverse can be used within GMRES-based iterative refinement (in place of the usual LU factorization) under certain constraints; here, we also have a constraint on how the working precision in iterative refinement (Algorithm \ref{alg:ir}) must be chosen relative to $\tau$. 

\begin{figure}[h]
\centering
\includegraphics[scale=.67]{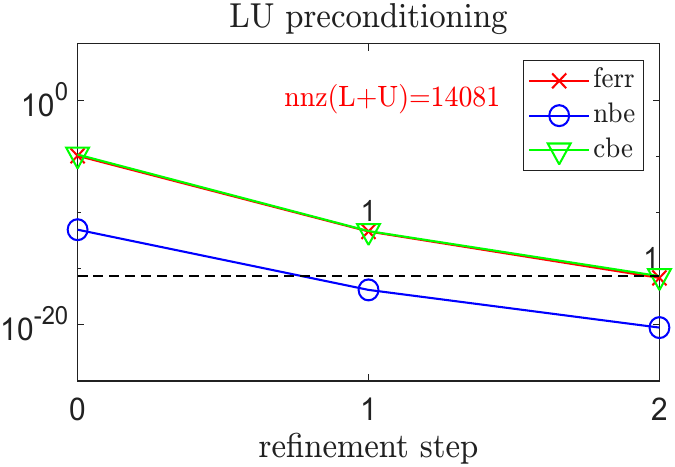}
\includegraphics[scale=.66]{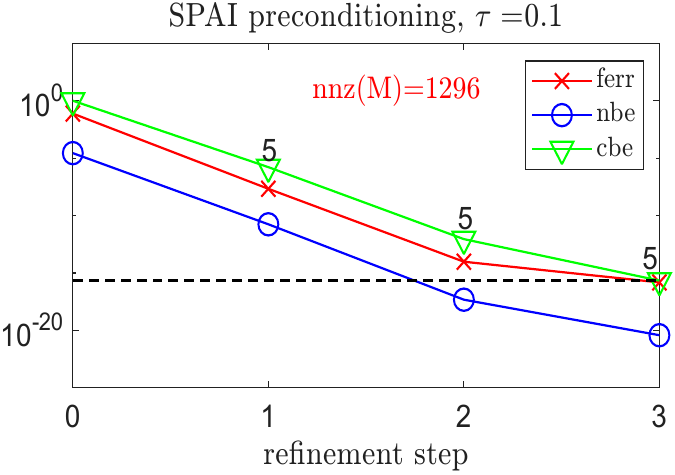}\\
\caption{Comparison of GMRES-IR with LU preconditioning (left) and SPAI preconditioning (right) for the matrix \texttt{steam1} from SuiteSparse \cite{davis2011university}. Numbers above markers give the number of GMRES iterations performed in each refinement step; in all cases the GMRES convergence tolerance was set to $10^{-6}$. For the top left plot, AMD ordering was applied before computing the LU factors. In all cases, we use GMRES-IR with $u_f=$ single, $u =$ double, $u_r=$ double. }
 \label{fig:SPAI}       
\end{figure}

In Figure \ref{fig:SPAI} we give one example of the use of low precision sparse approximate inverses within GMRES-based IR.  Here we use the matrix $\texttt{steam1}$ from SuiteSparse \cite{davis2011university}. We compare GMRES-based IR with LU preconditioning (left) and sparse approximate inverse preconditioning (right). The preconditioners (both LU and the sparse approximate inverse) were computed in single precision, and for the refinement process we use $u=u_r=$ double. In all cases, the GMRES convergence tolerance was set to $10^{-6}$. For SPAI, we set $\tau=0.1$. The plots in Figure \ref{fig:SPAI} show the forward error (`ferr'), normwise backward error (`nbe'), and componentwise backward error (`cbe'); the numbers above the markers indicate how many GMRES iterations were performed in each refinement step (i.e., in the solve in line \ref{line4} of Algorithm \ref{alg:ir}). In the LU and sparse approximate inverse plots, we list the size of the preconditioner in terms of number of nonzeros. We can see here that low-precision sparse approximate inverses are a viable alternative to LU-based preconditioners; while they require slightly more total GMRES iterations to converge, they result in a much less expensive preconditioner to apply and store. For further details and examples, see \cite{carson2023mixed}.

\subsection{Example 2: Mixed Precision Randomized Nystr\"{o}m Approximation}

Our next example involves computing a rank $k$ approximation of an $n\times n$ symmetric positive semidefinite matrix $A$ using a randomized Nystr\"{o}m approach. The randomized Nystr\"{o}m approximation has the form $A_N = (A\Omega)(\Omega^T A \Omega)^\dagger (A\Omega)^T$, where $\Omega$ is an $n\times k$ sampling matrix, and $\dagger$ indicates the Moore-Penrose pseudoinverse. The Nystr\"{o}m approximation arises in many applications, such as approximating kernel matrices and constructing spectral limited memory preconditioners. In some applications, the matrix is prohibitively large, and in some applications, its entries can only be accessed once. This motivated the single-pass variant of the Nystr\"{o}m method \cite{tropp2017fixed}, summarized in Algorithm \ref{alg:nys}.

\begin{algorithm}
\caption{Single-Pass Nystr\"{o}m Approximation \cite{tropp2017fixed}}\label{alg:nys}
\begin{algorithmic}[1]
\Require Symmetric positive semidefinite $A\in\mathbb{R}^{n\times n}$, target rank $k$
\Ensure $U\in\mathbb{R}^{n\times k}$, whose columns give approximate eigenvectors, $\Theta \in \mathbb{R}^{k \times k}$, whose diagonal entries give approximate eigenvalues
\State $G=\text{randn}(n,k)$, $[\Omega,\sim]=\text{qr}(G,0)$
\State $Y=A\Omega$ \label{nysline}
\State Compute shift $\nu$ and compute $Y_\nu = Y+\nu \Omega$, $B=\Omega^T Y_\nu$
\State $C=\text{chol}(B+B^T)/2$,  Solve $F=Y_\nu / C$
\State $[U,\Sigma,\sim]=\text{svd}(F,0)$, $\Theta = \text{max}(0,\Sigma^2-\nu I)$
\end{algorithmic}
\end{algorithm}

The matrix-matrix product in line \ref{nysline} is overwhelmingly the dominant cost in this regime. Our goal will thus be to use a working precision $u$ for all computations except this line, which is performed in precision $u_p\geq u$. Note that since this is the only time $A$ is accessed, $A$ should thus also be stored in precision $u_p$. Letting $\hat{A}_N$ denote the Nystr\"{o}m approximation computed in mixed precision, we can write 
$\Vert A-\hat{A}_N\Vert_2 \leq \Vert A-A_N \Vert_2 + \Vert A_N-\hat{A}_N\Vert_2$;
in other words, we bound the overall error as the sum of the exact approximation error and the finite precision error. Bounds for the first term have been previously derived; see, e.g., \cite{gittens2016revisiting}, \cite{frangella2023randomized}. In \cite{carson2024single}, we prove the first bound on the finite precision error. The result essentially says that as long as we choose precisions $u_p$ and $u$ such that $\kappa_2(A_k)\tilde{\kappa}(\Omega)^2 \ll u_p^{-2}$ and $\kappa_2(A_k)\tilde{\kappa}(\Omega)^2 \ll u^{-1}$, then 
\[
\Vert A_N - \hat{A}_N \Vert_F \lesssim \Vert A-A_N\Vert_F + k^{1/2} n u_p \kappa_2(A_k)\tilde{\kappa}(\Omega)^2\Vert A\Vert_F,
\]
where $A_k$ is the best rank-$k$ approximation of $A$ and $\tilde{\kappa}(\Omega)=\Vert \Omega \Vert_F \Vert (W_1^T \Omega)^\dagger)\Vert_2$, where $W_1$ are the eigenvectors for the leading $k$ eigenvalues of $A$. Under some assumptions we can take this bound and derive a heuristic which says that it is likely that the finite precision error is less than the approximation error when 
\[
u_p \leq n^{-1/2} \lambda_k / \lambda_1,
\]
where $\lambda_1$ and $\lambda_k$ are the largest and $k$th largest eigenvalues of $A$, respectively. 
This tells us again how the errors should be balanced: the greater the approximation error (the closer $\lambda_k$ is to $\lambda_1$), the lower the precision we can safely use. In Figure 2, we give one example of this for the matrix $\texttt{bcsstm07}$ from SuiteSparse \cite{davis2011university}. From the left plot, the heuristic above indicates that for $k\leq 148$, we can safely use half precision, and we can always safely use single precision. Indeed, this seems to be confirmed in the plot of the total error on the right; we see no difference with the exact case for single precision, and only start to see the finite precision error dominate for half precision around $k>150$. For more examples, an alternative (potentially more rigorous) heuristic, and insight on choosing $\Omega$ and $\nu$, see \cite{carson2024single}.

\begin{figure}[h]
\centering
\includegraphics[scale=.204]{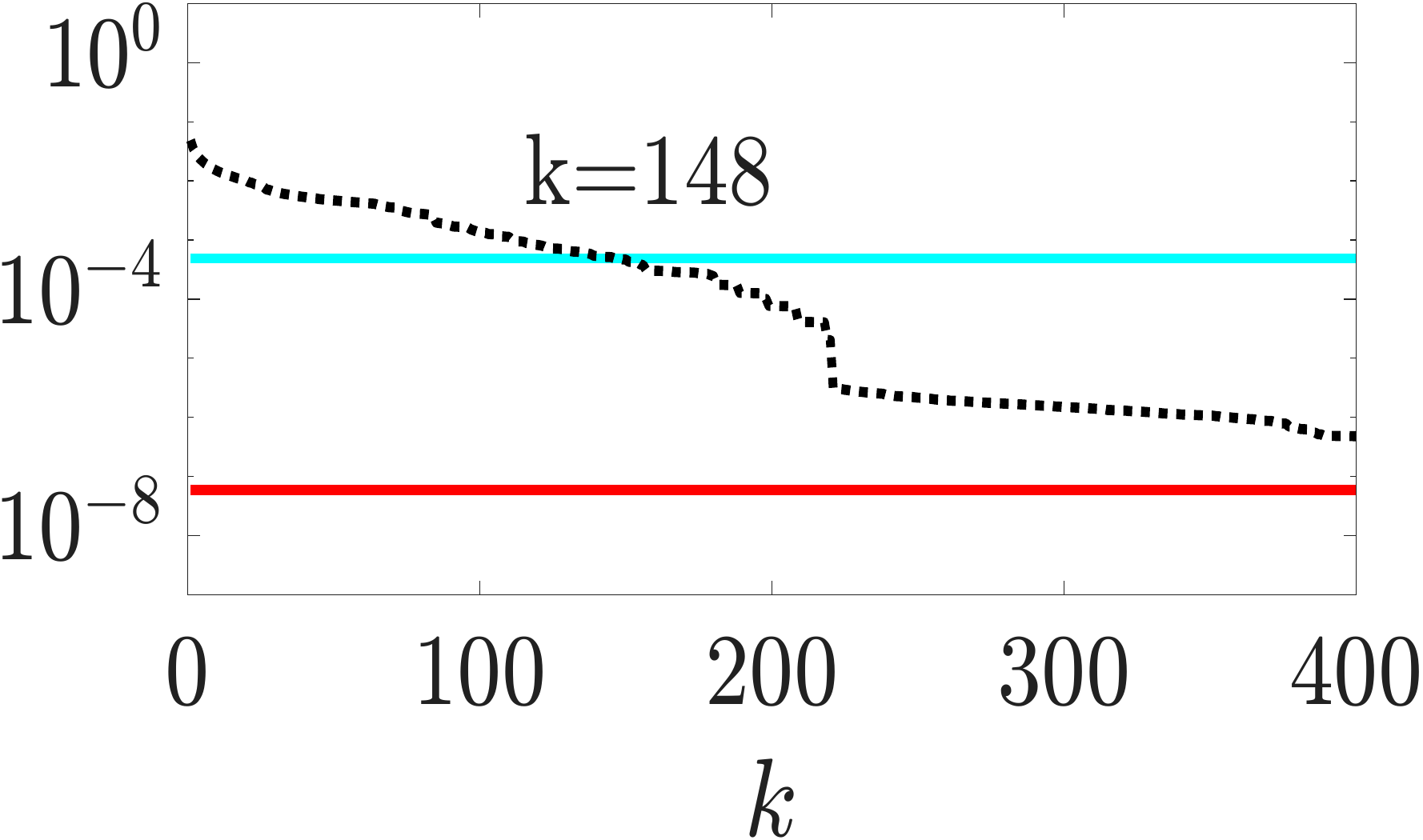}\hspace{4pt}
\includegraphics[scale=.2]{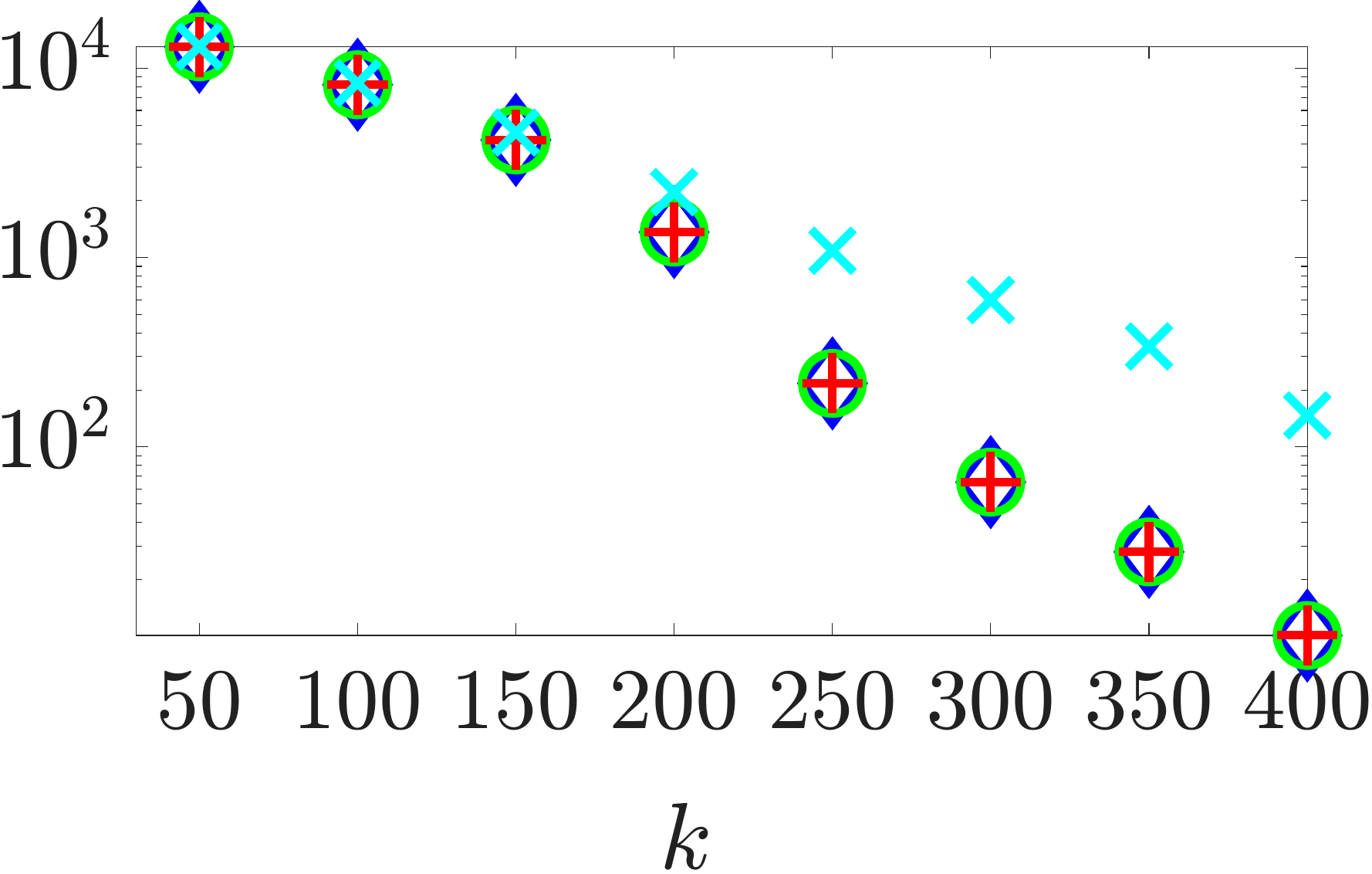}
\caption{Experiments for matrix \texttt{bcsstm07}. The left plot shows the scaled spectrum $\lambda_k/\lambda_1$ (black dotted), and the values of $\sqrt{n}u_p$ for half precision (cyan) and single precision (red). The heuristic indicates that the finite precision error should not dominate for $k$ less than the point where the black dotted line intersects the colored lines. The right plot shows the Frobenius norm of the mean total error, $\Vert A-\widehat{A}_N\Vert_F$, taken over five runs, for exact arithmetic (blue diamonds), $u_p,u=$ double (green circles), $u_p=$ single and $u=$ double (red $+$s), and $u_p=$ half and $u=$ double (cyan $\times$s). }
 \label{fig:Nystrom}      
\end{figure}

\subsection{Example 3: Mixed Precision Hierarchical Matrix Approximations} %

\begin{figure}[t]
\centering
{\includegraphics[trim= 30 350 30 0,clip,width=1\textwidth]{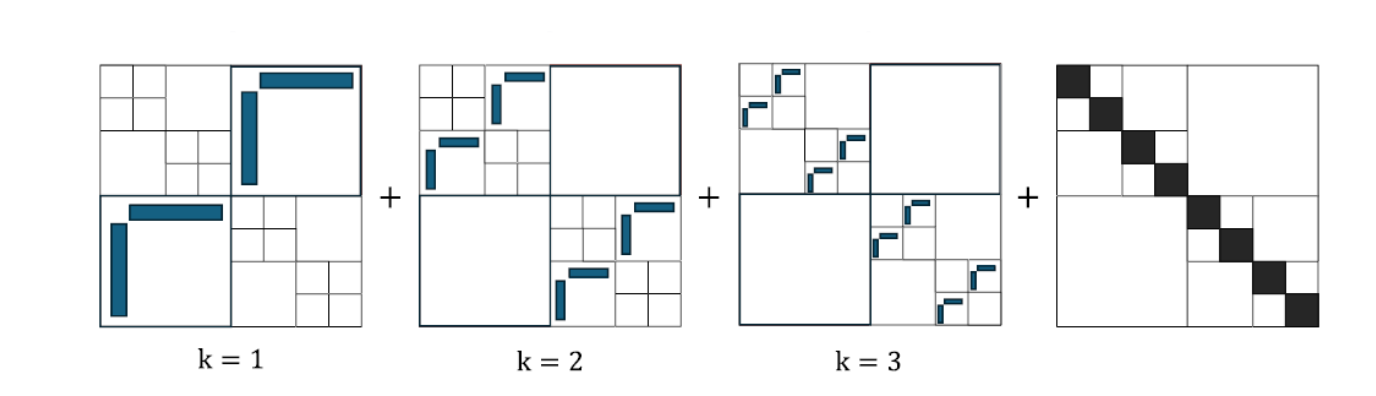}}
\caption{A depiction of a HODLR matrix with $\ell=3$ levels}
\label{fig:hodlr}
\end{figure}

The final example involves hierarchical matrix approximations, and in particular, HODLR (hierarchical off-diagonal low-rank) matrices. HODLR matrices have a fixed hierarchical block structure; given a number of levels $\ell$, we store the off-diagonal blocks in each level as rank-$p$ approximations, and the diagonal blocks in the final level as dense matrices. See the depiction with $\ell=3$ in Figure \ref{fig:hodlr}. Such representations reduce the cost of computations, e.g., $O(n^2)$ computations become $O(pn\log n)$.

The idea is that since we are already approximating the off-diagonal blocks, we can likely store them in lower precision as well without losing a significant amount of additional information. This is analyzed in \cite{carson2025mixed}.

Let $\varepsilon$ be the maximum relative error in the approximation of the off-diagonal blocks of a HODLR matrix $H$ (coming, for example, from using a truncated SVD).  
Let $\hat{H}$ be the mixed precision representation of $H$ where the off-diagonal blocks in level $k$ are stored in precision $u_k$ and the diagonal blocks in the final level are stored in a working precision $u<u_k$. Then we can write a bound for the global approximation error as 
\[
\frac{\Vert H - \hat{H}\Vert_F}{\Vert H\Vert_F} \lesssim 2\sqrt{2} \left( \sum_{k=1}^\ell 2^k \xi_k^2 u_k^2 \right)^{1/2} + \varepsilon,
\]
where 
$\xi_k = \max_{|i-j|=1} \Vert H_{ij}^{(k)} \Vert_F / \Vert H\Vert_F$ and $H_{ij}^{(k)}$, $|i-j|=1$ denotes any off-diagonal block from the $k$th level.

\begin{figure}[h]
\centering
\begin{subfigure}{0.2\textwidth}
    \centering
        \includegraphics[width=1.1\linewidth]{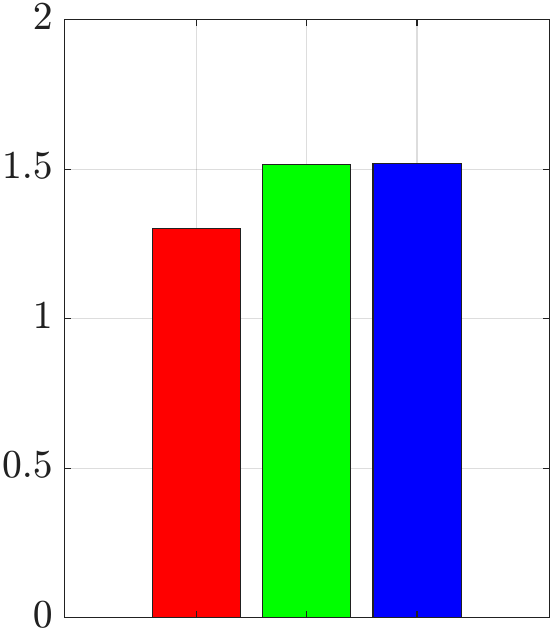}
        \caption{nos7}
        \label{fig:nos7}
\end{subfigure}%
\hspace{10pt}
\begin{subfigure}{0.2\textwidth}
    \centering
        \includegraphics[width=1.1\linewidth]{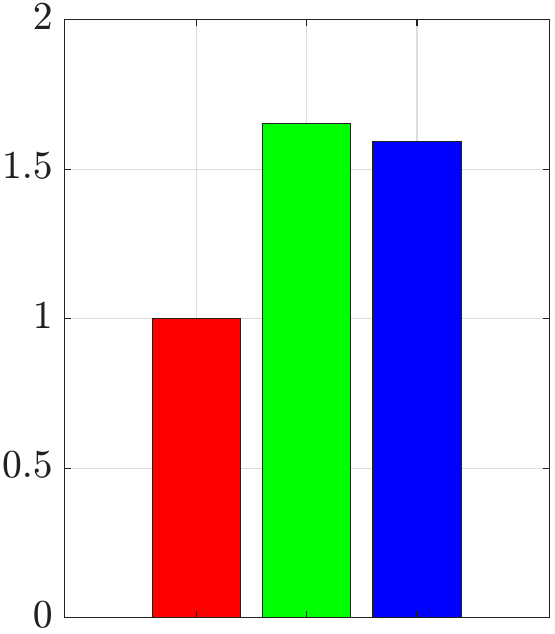}
        \caption{psmigr\_1}
        \label{fig:psm}
\end{subfigure}%
\hspace{10pt}
\begin{subfigure}{0.2\textwidth}
    \centering
        \includegraphics[width=1.1\linewidth]{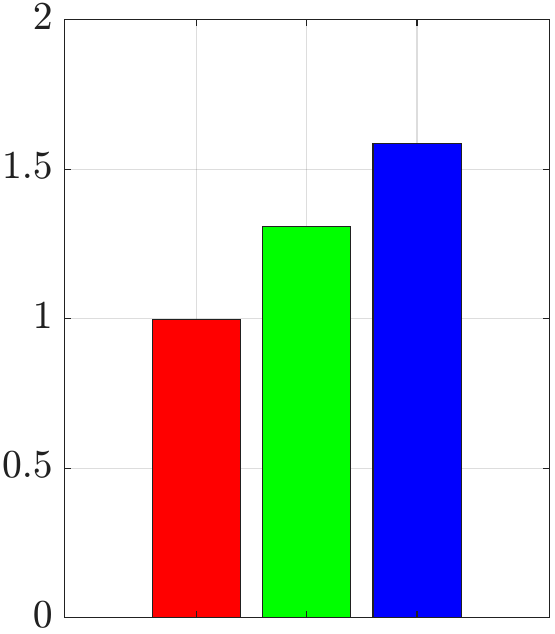}
        \caption{saylr3}
        \label{fig:saylr3}
\end{subfigure}%
\hspace{10pt}
\begin{subfigure}{0.2\textwidth}
    \centering
        \includegraphics[width=1.1\linewidth]{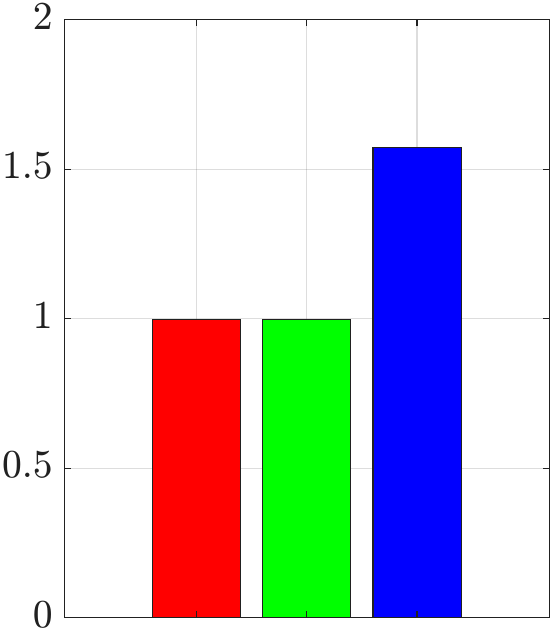}
        \caption{1138\_bus}
        \label{fig:1138bus}
\end{subfigure}%
\caption{Storage savings of adaptive-precision HODLR matrices relative to uniform (double) precision HODLR matrices. The depth $\ell=8$; red bars correspond to $\varepsilon=10^{-7}$, green bars indicate $\varepsilon=10^{-4}$, and blue bars indicate $\varepsilon=10^{-1}$.}
\label{fig:storage}
\end{figure}

This again tells us how to balance the errors: if we choose 
\[
u_k \leq \varepsilon/(2^{k/2}\xi_k),
\]
then the right-hand side of the above bound simplifies to $(2\sqrt{2} \ell +1)\varepsilon$; in other words, we do not see the effect of storing the off-diagonal blocks in lower precision. In \cite{carson2025mixed}, an adaptive precision approach to HODLR matrix construction is developed based on this idea. We show an example of the resulting storage savings for various $\varepsilon$ in Figure \ref{fig:storage} for a few problems from SuiteSparse \cite{davis2011university}.

In \cite{carson2025mixed} we also determine how to set the working precision $u$ in subsequent computations with mixed precision HODLR matrices so that the backward error of the computation does not exceed the error resulting from inexact representation of the matrix. The theorems say that for both matrix-vector products and LU factorizations with HODLR matrices stored in the mixed precision format, as long as $u \leq \varepsilon /n$, the backward error in these computations is dominated by the matrix representation error. Again, this gives us a theoretical guideline by which we can safely balance these errors. 
See \cite{carson2025mixed} for details and examples.

\section{Other Examples and Outlook}
\label{sec:3}

There are many other examples of balancing finite precision error with other errors, e.g., discretization and sampling errors \cite{martinek2025exploiting}, and hardware failures \cite{carson2025detection}, that we did not have space to include here. 
For readers interested in mixed precision numerical algorithms, we refer to the surveys \cite{abdelfattah2021survey} and \cite{higham2022mixed}.

Looking forward, it is clear that we are moving more and more towards the era of extreme heterogeneity in large-scale systems, meaning that it is more important than ever to take a holistic approach to algorithm design, and truly consider the whole computational science process, including the hardware we run our codes on. 

There is a great amount of work going on in developing new non-IEEE number formats and arithmetic systems, developing approximate hardware, and the emergence of new computing paradigms like quantum computing, neuromorphic computing and their hybrids. The task of analyzing errors from multiple sources and determining how to balance errors in an optimal way in these new paradigms will surely provide a source of interesting, challenging problems for a long time to come. 

We end with a challenge to the reader to consider their own application area, the sources of error involved, and where it might make sense to use low or mixed precision computations.

\paragraph{Acknowledgements.} This manuscript is dedicated to the memory of Nicholas J. Higham, a mentor and friend, who first piqued my interest in mixed precision algorithms. 

\bibliographystyle{plain}
\bibliography{arxivversion.bbl}

\end{document}